\documentclass[11pt]{article}
\usepackage{amsxtra,amssymb,amsfonts,amsthm,amsmath}
\usepackage{graphicx}
\usepackage{diagrams}

\theoremstyle{plain}
\newtheorem{theorem}{Theorem}
\newtheorem*{theorem*}{Theorem}
\newtheorem{lemma}{Lemma}

\newcommand{\refT}[1]{Theorem~\ref{T:#1}}
\newcommand{\refS}[1]{Section~\ref{#1}}
\newcommand{\refL}[1]{Lemma~\ref{L:#1}}

\newcommand{\refF}[1]{Fig.~\ref{F:#1}}

\def\ds{\displaystyle}

\def\R{{\mathbb R}}

\def\eps{\varepsilon}

\begin{document}

\title{ {\LARGE {\bf Invasion waves in the presence of a mutualist} }}

\author{
{\sc Volodymyr Hrynkiv and Sergiy Koshkin } \\
Computer and Mathematical Sciences Department\\
University of Houston Downtown\\
Houston, TX 77002\\
}

\date{}

\maketitle

\begin{abstract}

This paper studies invasion waves in the diffusive
Competitor-Competitor-Mutualist model generalizing the $2$-species
Lotka-Volterra model studied by Weinberger et al. The mutualist
may benefit the invading or the resident species producing two
different types of invasions. Sufficient conditions for linear
determinacy are derived in both cases, and when they hold,
explicit formulas for linear spreading speeds of the invasions are
obtained by linearizing the model. While in the first case the
linear speed is increased by the mutualist, it is unaffected
in the second case. Mathematical methods are based on
converting the model into a cooperative reaction-diffusion system.
\end{abstract}

\vspace{.12in} \noindent {\it Keywords}: cooperative system, reaction-diffusion, spreading speed, monostable waves, linear determinacy, Competitor-Competitor-Mutualist, Lotka-Volterra.\\


\section*{Introduction}\label{s0}

Traveling waves appear in many biological models with spatial structure. The first such model, for a spread of an advantageous gene, was
introduced by Fisher in 1937. It leads to the Fisher equation $\dot{p}=rp(1-p)+Dp_{xx}$, where $p$ is the fraction of the local population that contains the advantageous gene, $r$ is its initial growth rate, and $D$ measures its mobility. This equation has two constant equilibria, an unstable one at $p=0$
and a stable one at $p=1$, and admits a traveling wave solution $p=w(x-ct)$ with the profile $w$ approaching $1$ at $-\infty$ and $0$ at $\infty$. The wave speed $c=2\sqrt{rD}$ is also the asymptotic spreading speed for this equation. Namely, for any $p(0,x)$ with compact support and $\eps>0$ we have
$$\ds{\lim_{t\to\infty}p(t,x)=\begin{cases}0,&|x|\geq(c+\eps)t\\
1,&|x|\leq(c-\eps) t.
\end{cases}}$$

More recently, traveling waves were studied in multi-species cooperative and competition models. For example, the
$2$-species diffusive Lotka-Volterra competition model
\begin{equation}\label{2LotVl}
\begin{cases}\dot{p}_1=\alpha p_1(1-p_1-ap_2)+D_1p_{1,xx}\\
\dot{p}_2=\beta p_2(1-p_2-bp_1)+D_2p_{2,xx}\end{cases}
\end{equation}
is analyzed in \cite{W2'}. Here $p_i$ are the population densities of two competing species, and the constants 
$a,b$ characterize the intensity of the competition. In \cite{W2'} the first species is interpreted as
a non-native invader of a habitat, originally occupied exclusively by the second, resident species. It is shown that in a certain range of
parameters there is again a traveling wave connecting the invaded equilibrium $(0,1)$ either to $(1,0)$, where the resident
is extinct, or to a coexistence equilibrium $(p_1^*,p_2^*)$. The asymptotic spreading speed is explicitly found to be $c=2\sqrt{\alpha D_1(1-a)}$
under conditions that guarantee the so-called linear determinacy. {\bf Linear determinacy} means that the non-linear system \eqref{2LotVl} has the same
spreading speed as its linearization at the invaded equilibrium $(0,1)$, i.e. as the linearization at the edge of the spreading wavefront, provided that
the initial density vector is equal to $(0,1)$
outside a bounded subset.
The analysis of the system \eqref{2LotVl} in \cite{W2'} uses the substitution $q_1=p_1$, $q_2=1-p_2$, which converts it into a cooperative parabolic system
of reaction-diffusion type.

Recall that a system $\dot{p}_i=\beta
f(p_1,\dots,p_n)+D_ip_{i,xx}$ is cooperative if $\frac{\partial
f_i}{\partial p_j}\geq0$ for $i\neq j$, i.e. small increase in any
species' population is beneficial to all other species. Weinberger
et al. showed in \cite{W2} that analysis of such a system can be
reduced to analysis of recursions $p^{(n+1)}=Q[p^{(n)}]$, where
$Q$ is the time $1$ map of the evolution system, and gave
sufficient conditions for linear determinacy to hold (later
revised in \cite{W5,W7}).

The progress of an invasion by a competitor may depend on many additional
factors not captured by a $2$-species model. For instance,
the interaction of both the invader and the resident with a third
species may play a role if its presence is beneficial to one or
both of them. In this paper we examine the influence of such
species, a mutualist, on invasion waves. As a model we pick the
Competitor-Competitor-Mutualist model introduced in \cite{RFA}
with diffusion added:
\begin{equation}\label{CCM}
\begin{cases}\dot{p}_1=\alpha p_1(1-p_1-\frac{ap_2}{1+mu})+D_1\,p_{1,xx}\\
\dot{p}_2=\beta p_2(1-p_2-bp_1)+D_2\,p_{2,xx}\\
\dot{u}\,\,=\gamma u(1-\frac{u}{L+lp_1})+D_3\,u_{xx}.
\end{cases}
\end{equation}
The population density of the mutualist is $u$, whereas $m$ and $l$
measure intensities of mutualism for the first species and the
mutualist, respectively. In this model the mutualist benefits the
first competitor by mitigating competition with the second one.
The mutualist's growth is logistic with the carrying capacity
$L+lp_1$, so $L$ is the self-carrying capacity in the absence of
the first species. We always assume $L>0$ because the last equation is singular at the origin otherwise. 
In the absence of competition ($a=0$) there is
no cost or benefit to the first competitor, and there is no direct
interaction between the mutualist and the second competitor. Also
note that if $m=0$ the first two equations decouple and form a
Lotka-Volterra competition system \eqref{2LotVl}.

We consider two situations. In the first one, the mutualist
benefits the invader, i.e. $p_1$ describes the invader and $p_2$
the resident. For instance, $u$ may describe a microorganism that
helps the invader compete with the resident species. One example is 
the role of mycorrhizal fungi in the competition between saltcedar 
and cottonwood \cite{Beau, Beau1}. Such a {\bf mutualist-invader}, as we shall call it, can also be used as a
biocontrol agent if the invasion is desirable, e.g. if the
resident is considered a pest \cite{L-T}. In the second situation, a
{\bf mutualist-resident} mitigates an invasion by a non-native
species by reducing the competition with the invader, i.e. $p_1$
describes the resident and $p_2$ the invader. Such a mutualist may
be used to stem an invasion by pests \cite{BD}. In this paper we do not address an interesting third possibility of a mutualist invading a coexistence equilibrium of the first two species.

In each case we consider,
appropriate equilibria are selected as initial values, and
sufficient conditions for the existence and linear determinacy of
invasion waves are derived. We also estimate the spreading
speeds and compare the results to those for the purely competitive
model \eqref{2LotVl}. As in \cite{W2'}, our approach is
based on performing substitutions that convert \eqref{CCM} into
cooperative systems, and applying to them the analytic apparatus
of \cite{W2}. Our main results are summarized in Theorems
\ref{T:invader} and \ref{T:resident}.

The paper is organized as follows. \refS{s3} is a brief introduction to the theory of spreading
speeds in cooperative reaction-diffusion systems developed by Lui,
Weinberger and others. In \refS{s1} we recall and
refine the results of \cite{RFA} on equilibria of the
Competitor-Competitor-Mutualist system without diffusion. They
serve as initial and asymptotic states for invasion waves in
system \eqref{CCM}. Sections \ref{s4} and \ref{s5} present our
analysis of invasion waves in system \eqref{CCM} for the
mutualist-invader and the mutualist-resident, respectively. 
Our simulations that estimate the spreading speeds numerically are described in \refS{s6}. 
We conclude with a discussion in \refS{s7}. Proofs of technical lemmas are collected in the Appendix.

\section{Spreading speeds in cooperative systems}\label{s3}

Our approach to analyzing the invasion waves in the diffusive
Competitor-Competitor-Mutualist system is based on converting it
into cooperative systems, for which the theory is well-developed.
The conversion is performed in Sections \ref{s4} and \ref{s5} for the mutualist-invader and the mutualist resident
case, respectively. In this section we give a brief overview of the theory of the spreading speeds in cooperative systems
due to Weinberger et al. \cite{W2}, as it applies to reaction-diffusion systems.

Let us write $p\leq
q$ and $p\ll q$ for vectors $p$, $q$ if $p_i\leq q_i$ and
$p_i<q_i$ respectively, for all $i$. The same notation is used for
vector functions if the same inequalities hold pointwise. 
Denote $\mathcal{V}_{\beta}:=\{p\in\R^n|\ 0\leq
p\leq\beta\}$ and $\mathcal{C}_{\beta}$ the space of continuous
functions with values in $\mathcal{V}_{\beta}$. Consider a reaction-diffusion equation 
for a vector function $p(t,x)$:
\begin{equation}\label{RD}
\dot{p}=f(p)+Dp_{xx}.
\end{equation}
The system \eqref{RD} is called cooperative on
$\mathcal{V}_{\beta}$ if in this region $\frac{\partial
f_i}{\partial p_j}\geq0$ for $i\neq j$, i.e. the Jacobian matrix
has positive entries off-diagonal \cite{Sm}. Let $Q_t$ be the
evolution operator for \eqref{RD}, i.e. $Q_t[p](x):=p(t,x)$ is
the solution to the system with $p(0,x)=p(x)$. Cooperativity
implies that the evolution is monotone: within
$\mathcal{V}_{\beta}$ the inequality $p(x)\leq q(x)$ implies
$Q_t[p](x)\leq Q_t[q](x)$ for all $t>0$. In particular,
$\mathcal{C}_{\beta}$ is invariant under $Q_t$ since the
inequalities $0\leq p(x)\leq\beta$ will be preserved by it.

Suppose $f(0)=0$ so that the linearization of \eqref{RD} at $p=0$ is
\begin{equation}\label{LinRD}
\dot{p}=f'(0)p+Dp_{xx}.
\end{equation}
Even for the linearized equation different components of $p$ may
spread at different speeds \cite{ILN,W7}. If $f'(0)$ has a strictly
positive eigenvector with a positive eigenvalue, and if there are
no strictly positive equilibria between $0$ and $\beta$ standard
monotonicity arguments show that any solution to $\dot{p}=f(p)$
with $p(0)=x$ satisfying $0\ll x\leq\beta$ converges to $\beta$.
In other words, the equilibrium $\beta$ is stable in the large for
solutions initiating in the interior of
$\mathcal{V}_{\beta}:=\{p\in\R^n|\ 0\leq p\leq\beta\}$. Moreover,
it is shown in \cite{W5,W2} that for solutions to \eqref{RD} with
$p(0,x)$ having compact support and strictly positive on a
sufficiently long interval there exist the slowest and the
fastest spreading speeds. The {\bf slowest spreading speed $c^*$} is characterized by requiring 
that for every $\eps>0$:
\begin{equation}\label{c*}
\lim_{t\to\infty}\sup_{|x|\leq(c*-\eps)t}\Big(\beta-p(t,x)\Big)=0\,,
\hspace{1em}\text{and for some $i$}\hspace{1em}
\lim_{t\to\infty}\sup_{|x|\geq(c*+\eps)t}p_i(t,x)=0\,.
\end{equation}
That is, the $i^{\text{th}}$ component spreads at a speed no higher than $c^*$, and no other component spreads at lower speed. Analogously, the {\bf fastest spreading speed $c^*_f$} is characterized by
\begin{equation}\label{c*f}
\lim_{t\to\infty}\sup_{|x|\geq(c^*_f+\eps)t}p(t,x)=0\,,
\hspace{1em}\text{and for some $i$}\hspace{1em}
\limsup_{t\to\infty}\inf_{|x|\leq(c^*_f-\eps)t}p_i(t,x)>0\,,
\end{equation}
i.e. the $i^{\text{th}}$ component spreads at a speed no less than $c^*_f$, and no component spreads at higher speed.
Clearly, $c^*\leq c^*_f$ and the system is said to have a {\bf single speed} if $c^*=c^*_f$.
A simple but useful consequence of \eqref{c*}, \eqref{c*f} is that if two solutions $p,\widetilde{p}$ satisfy $p(t,x)\leq\widetilde{p}(t,x)$ and both have spreading speeds, then those are related by the same inequality: 
$c^*\leq\widetilde{c}^*$ and $c^*_f\leq\widetilde{c}^*_f$.

In the case of the linearized system \eqref{LinRD} the spreading speeds
can be explicitly computed and used to estimate the non-linear
speeds. To proceed we need some more terminology. By permuting
coordinates it is possible to bring $f'(0)$ to the {\bf Frobenius
form}, i.e. a lower block-triangular form with irreducible
diagonal blocks. The irreducibility means that permutations can not
block-triangulate those blocks any further. Consider the matrix 
\begin{equation}\label{Cmu}
C_\mu:=f'(0)+\mu^2D.
\end{equation}
\noindent Since $D$ is diagonal $C_\mu$ will be in the Frobenius form whenever $f'(0)$
is. By a theorem of Perron-Frobenius \cite{Th}, each irreducible
diagonal block has a single real eigenvalue with the eigenvector having strictly
positive entries within this block. It is called the {\bf principal
eigenvalue.} Let $\gamma_i(\mu)$ be the
principal eigenvalue of the $i^{th}$ block and $\zeta_i(\mu)$ be
the corresponding eigenvector of the entire matrix $C_\mu$ (it may not
have all positive entries in general). Then
$p(t,x)=e^{\gamma_i(\mu)t-\mu x}\zeta_i(\mu)$ is a traveling wave
solution to \eqref{LinRD}. One can now define the minimal wave speeds
for each block as $\ds{c_i:=\inf_{\mu>0}\mu^{-1}\gamma_i(\mu)}$, they
would be the spreading speeds for the block's components if all
other components are kept at $0$ \cite{Lui}. Even if these {\bf
block speeds} are different it is still possible that the
linearized system has a single speed. Consider the {\bf linear
speeds} $c_1$ and $\ds{c_f:=\max_ic_i}$. A system is called {\bf
linearly determinate} if $c_1=c^*$ and $c_f=c^*_f$. Interestingly
enough, even for a linear system it can happen that $c^*_f>c_f$,
see \cite{W7}. The following theorem follows from the results of
\cite{W5,W2}.

\begin{theorem}\label{T:LDet} Let $f$ be a function with $f(0)=0$, $f(\beta)=0$ for $\beta\gg0$, and no $\nu\gg0$  between $0$ and $\beta$.
Suppose that $f$ is smooth and cooperative on
$\mathcal{V}_{\beta}:=\{p\in\R^n|\ 0\leq p\leq\beta\}$. Let $D$ be a diagonal matrix with positive entries, and let $f'(0)$and hence $f'(0)+\mu^2D$ be
in the Frobenius form with the eigenvalues $\gamma_i(\mu)$ and the eigenvectors
$\zeta_i(\mu)$. If

{\rm I)} $\gamma_1(0)>0$ and has a strictly positive eigenvector $\zeta_1(0)\gg0$

\noindent then $c^*$, $c^*_f$ exist for $\dot{p}=f(p)+Dp_{xx}$ and $c^*\geq c_1 $, $c^*_f\geq c_f$. If also

{\rm II)} $\gamma_1(\overline{\mu})>\gamma_i(\overline{\mu})$ for
$i>1$, where $\overline{\mu}$ is the number at which
$\ds{\inf_{\mu>0}\mu^{-1}\gamma_1(\mu)}$ is attained;

{\rm III)} $f(\rho\zeta_1(\overline{\mu}))\leq\rho f'(0)\zeta_1(\overline{\mu})$ for all real numbers $\rho>0$;

\noindent then the system has a single spreading speed and is linearly determinate:\\
$\ds{c^*_f=c^*=c_f=c_1=\inf_{\mu>0}\mu^{-1}\gamma_1(\mu)}$.
\end{theorem}
\noindent Let us briefly comment on the nature of the assumptions
I)--III), see \cite{W2} for more details. Growth of invading
species at the zero density is usually unaffected by the presence of
other species, meaning that the off-diagonal entries in the
corresponding row of $f'(0)$ vanish. This implies that the
invader's diagonal block is $1\times1$, and it has to be placed first
in the Frobenius form. Our first assumption then means that there
is a single invader and at zero initial density its density grows.
Second assumption requires the minimal block speed of the invader
to dominate the minimal block speeds of all other species.
Finally, III) roughly means that the growth rates decrease with
population density (no Allee effect), at least along the direction
of the invader's eigenvector.

Note that the theorem only guarantees the existence of a spreading speed for initial densities
with compact support. There may or may not exist a traveling wave profile $w(x)$ like the one in the Fisher model,
i.e. with $w(-\infty)=\beta$, $w(\infty)=0$ and $w(x-ct)$ satisfying the equation. The subtlety is that although there are no strictly positive equilibria between $0$ and $\beta$ there may be boundary ones with one or more zero components. According to \cite{W5}, there will always be a wave profile with $w(-\infty)=\beta$ and $w(\infty)$ equal to such a {\bf boundary equilibrium}, but not necessarily to $0$. Boundary equilibria may also induce complicated dynamics such as development of stacked wavefronts propagating at different speeds \cite{ILN}. In their absence not only does Fisher-like profile exist, but \eqref{RD} also has a single spreading speed even if assumptions II)-III) are dropped. However, that speed may not be linearly determinate.

\section{Constant equilibria}\label{s1}

In this section we discuss the nature and the relative positions of the spatially constant equilibria in  the Competitor-Competitor-Mutualist system.
For convenience, we identify $p_3=u$ and use both notations interchangeably. Introducing the column vector $p:=\begin{pmatrix}p_1&p_2&u\end{pmatrix}^T$
and the diagonal matrix $D:=\text{diag}(D_i)$ we can rewrite the system \eqref{CCM} in the vector form $\dot{p}=f(p)+Dp_{xx}$ with
\begin{equation}\label{fCCM}
f(p)=\begin{pmatrix}\alpha\,p_1(1-p_1-\frac{ap_2}{1+mu})\\
\beta\,p_2(1-p_2-bp_1)\\
\gamma\,u(1-\frac{u}{L+lp_1})\end{pmatrix}.
\end{equation}
The constant equilibria are the solutions to $f(p)=0$. Listed in the order from \cite{RFA} they are
\begin{align}\label{Eqbr}
E_0\,(0,0,0) & & E_1\,(1,0,0) & & E_2\,(0,1,0) & & E_3\,(0,0,L) \notag\\
E_4\,(0,1,L) & & E_5\!\left(\frac{1-a}{1-ab},\frac{1-b}{1-ab},0\right) & & E_6\,(1,0,L+l)& & 
E_7^\pm\,(z_\pm,1-bz_\pm,L+lz_\pm)\,,
\end{align}
where $z_\pm$ are solutions to the quadratic equation
\begin{equation}\label{zEq}
ml\,\!z^2+(1-ab+mL-ml)z-(1+mL-a)=0.
\end{equation}
The listed {\bf equilibria are admissible}
only if their entries are non-negative since they represent the
population densities of the corresponding species. In particular, the
admissibility of the coexistence equilibria $E_5$ and $E_7^\pm$
depends on values of the parameters. As the notation suggests, $E_7^\pm$ may in fact
represent up to two admissible equilibria depending on the
number of solutions to \eqref{zEq} and their values.

A constant {\bf equilibrium is invadable} by a species if its population is zero at the equilibrium, but grows if the zero is replaced by a small positive value \cite{W2}. In most cases the system then evolves into another equilibrium, so
there is an orbit connecting the two equilibria of the ODE system $\dot{p}=f(p)$ with
$f(p)$ given by \eqref{fCCM}. We call them {\bf the source and the target equilibria} respectively, they are the states at which the orbit originates and which it asymptotically approaches. A traveling wave profile, if it exists as in the Fisher model, will approach the source equilibrium at $\infty$ and the target one at $-\infty$.
If more than one species is missing from an equilibrium $E_i$ then one computes from \eqref{fCCM} that $f'(E_i)$ has two 
$1\times1$ diagonal blocks violating a basic condition of \cite{W2} on the absence of multiple invaders. 

This leaves $E_4,E_5$ and $E_6$ as the possible source equilibria, where the methods of Weinberger et al. can be applied. It may seem that $E_7^\pm$ is also a possibility when $z_\pm=1/b$, but one can check that then $b=1$, and hence $E_7^\pm=E_6$. When $E_4$ is the source equilibrium the mutualist benefits the absent species, i.e. the invader, and we refer to it as a {\bf mutualist-invader}. Accordingly, a {\bf mutualist-resident} appears with
$E_6$ as the source since the mutualist then benefits the species present in the equilibrium. In contrast, when $E_5$ is the source equilibrium it is the mutualist that invades the equilibrium of two coexisting species, a situation with no analog in the $2$-species Lotka-Volterra competition model. We do not consider it in this paper.

\section{Mutualist-invader}\label{s4}

In this section we derive explicit sufficient conditions for a single speed and linear determinacy
from \refT{LDet} in the case of a mutualist-invader. This means that the source equilibrium is $E_4$, and first we perform a substitution that moves the source equilibrium to the origin and turn the system into a cooperative one. The target equilibrium then serves as $\beta$ from Section \ref{s3}. Following Weinberger et al., we will also assume that 
$\beta$ is strictly positive, and there are no other strictly positive equilibria between $0$ and $\beta$. Under these conditions an orbit connecting $0$ to $\beta$ is guaranteed to exist \cite{W2}.

For the mutualist-invader we transform \eqref{CCM} by extending the corresponding substitution for the $2$-species model from \cite{W2'}:
$$\ds{\left[\begin{array}{l} q_1=p_1\\
q_2=1-p_2\\
v=u-L\ (q_3=p_3-L).
\end{array}\right.}$$
As before we identify $v=q_3$ and set $q:=\begin{pmatrix}q_1&q_2&v\end{pmatrix}^T$\!\!, $D:=\text{diag}(D_i)$. The transformed system in the vector form is $\dot{q}=g(q)+Dq_{xx}$ with
\begin{equation}\label{gCCMI}
g(q)=\begin{pmatrix}\alpha\,q_1(1-q_1+\frac{a(q_2-1)}{1+mL+mv})\\
\beta\,(1-q_2)(bq_1-q_2)\\
\gamma\,(v+L)\left(1-\frac{v+L}{L+lq_1}\right)\end{pmatrix}.
\end{equation}
A straightforward computation for the Jacobi matrix $\ds{g'(q):=\left(\frac{\partial g_i}{\partial q_j}\right)}$ yields
\begin{equation}\label{g'CCMI}
g'(q)=\begin{bmatrix}
\alpha\left(1-2q_1+\frac{a(q_2-1)}{1+mL+mv}\right) & \alpha\,\frac{aq_1}{1+mL+mv} & \alpha\,\frac{maq_1(1-q_2)}{(1+mL+mv)^2}\\
\beta\,b(1-q_2)& \beta\,(2q_2-bq_1-1) & 0\\
\gamma\,l\frac{(v+L)^2}{(L+lq_1)^2} & 0 & \gamma\left(1-2\frac{v+L}{L+lq_1}\right)\end{bmatrix}.
\end{equation}
One can see by inspection that for non-negative values of $\alpha,\beta,\gamma,a,b,m,L,l$ and $q_2\leq1$ the off-diagonal entries in \eqref{g'CCMI} are also non-negative, i.e. the transformed system is
cooperative. 

The constant equilibria are transformed into (cf. \eqref{Eqbr}):
\begin{align}\label{Fqbr}
F_0\,(0,0,-L) & & F_1\,(1,1,-L) & & F_2\,(0,0,-L) & & F_3\,(0,1,0)\notag\\
F_4\,(0,0,0)& & F_5\!\left(\frac{1-a}{1-ab},b\frac{1-a}{1-ab},-L\right) & & F_6\,(1,1,l)& & F_7^\pm\,(z_\pm,bz_\pm,lz_\pm),
\end{align}
as before $z_\pm$ is a non-negative solution to \eqref{zEq}. The source equilibrium $E_4$ moves to the
origin $F_4$. By inspection, the only strictly positive target equilibria are
$F_6$ and possibly $F_7^{\pm}$. For the original system they correspond to 
the resident extinction $E_6\,(1,0,L+l)$, and the $3$-species
coexistence $E_7^+\,(z_+,1-bz_+,L+lz_+)$, respectively. 

Since \refT{LDet} mostly involves conditions on $g'(0)$ we compute
from \eqref{g'CCMI}
\begin{equation}\label{g'0CCMI}
g'(0)=\begin{bmatrix}
\alpha\left(1-\frac{a}{1+mL}\right) & 0 & 0\\
\beta\,b & -\beta & 0\\
\gamma\,l & 0 & -\gamma
\end{bmatrix}.
\end{equation}
This matrix is already in the Frobenius form with $1\times1$ diagonal blocks. For the minimal block speeds to exist at all we need at least one positive eigenvalue.
It can only be $\gamma_1(0)=\alpha\left(1-\frac{a}{1+mL}\right)$ since the initial growth rates $\beta$, $\gamma$ are always positive, and it is positive only if $a<1+mL$. The following lemma provides additional information on the coexistence equilibria $E_7^\pm$ when this last inequality is satisfied. \refL{a1mL} is essentially proved in \cite{RFA}, but for the convenience of the reader we give another proof in the Appendix.
\begin{lemma}\label{L:a1mL} Suppose $a,m,l,b,L>0$ and $a<1+mL$. Then $E_7^-$ is never admissible and $E_7^+$ is admissible if and only if $b\leq1$. If $b=1$ then $E_7^+=E_6$, otherwise $E_7^+$ is distinct from the other equilibria.
\end{lemma}
\noindent It follows from \refL{a1mL} that $F_7^\pm$ are inadmissible for $b>1$, and
from its proof it follows that $F_7^+\ll F_6$ for $b<1$ (because then also $z_+<1$, see \eqref{1-z+}). We conclude that for $b<1$ the target is $F_7^+$ and for $b>1$
it is $F_6$, for $b=1$ they are equal, $F_6=F_7^+$. In both cases the second component of the target equilibrium is 
$\leq1$. Since $\mathcal{C}_\beta$ with $\beta=F_6,\,F_7^+$ is invariant under the time evolution we conclude that the value of $q_2$ will remain $\leq1$ for the duration, keeping the system cooperative.

Even without linear determinacy we can draw some conclusions on a single speed and the existence of traveling waves. As explained above, for $b<1$ the target equilibrium is $F_7^+$, which corresponds to the $3$-species coexistence $E_7$ in the original system. Moreover, $z<1$ in $F_7^+$
and there are no boundary equilibria in
$\mathcal{C}_{\beta}$. It follows from \cite{W5} that there is a traveling wave profile $w$ with $w(-\infty)=\beta$
and $w(\infty)=0$, and the system has a single asymptotic spreading speed equal to the wave speed of this profile. 
If $b<1$ then there is a
boundary equilibrium $F_5$. For $b\geq1$ the target equilibrium switches to $F_6$, corresponding to the resident extinction, and a boundary equilibrium
$F_3$ is acquired. In the latter cases we can not make definitive conclusions on a single speed without linear determinacy.

Going back to linear determinacy, the matrix $C_\mu:=g'(0)+\mu^2D$ from \eqref{Cmu} is obtained by adding $D_i\mu^2$ to the diagonal entries of $g'(0)$, and the obtained entries $\gamma_i(\mu)$ are its eigenvalues:
$$
\gamma_1(\mu)=\alpha\left(1-\frac{a}{1+mL}\right)+D_1\mu^2;\hspace{2em} \gamma_2(\mu)=-\beta+D_2\mu^2;\hspace{2em}
\gamma_3(\mu)=-\gamma+D_3\mu^2.
$$
The principal eigenvalue is $\gamma_1$ and the corresponding eigenvector is
\begin{equation}\label{z1(0)}
\zeta_1=\begin{pmatrix}
(\gamma_1-\gamma_2)(\gamma_1-\gamma_3)\\
\beta\,b\,(\gamma_1-\gamma_3)\\
\gamma\,l\,(\gamma_1-\gamma_2)
\end{pmatrix},
\end{equation}
where we suppressed from notation the dependence of both sides on $\mu$. One can see by inspection that
$\zeta_1(0)$ has positive entries as long as $a<1+mL$.

Assuming the system is linearly determinate the linear speed will be
\begin{equation}\label{c1i}
\ds{c_1=\inf_{\mu>0}\mu^{-1}\gamma_1(\mu)=2\sqrt{\alpha D_1\left(1-\frac{a}{1+mL}\right)}}
\end{equation}
with the infimum attained at
\begin{equation}\label{mubar}
\overline{\mu}=\sqrt{\frac{\alpha}{D_1}\left(1-\frac{a}{1+mL}\right)}.
\end{equation}
Note that for the self-carrying capacity of the mutualist $L$ near zero the speed is nearly the same as in 
the $2$-species model \eqref{2LotVl}. As $L$ increases the speed goes up to $\ds{c=2\sqrt{\alpha D_1}}$, where the effects of the competition are completely erased by the influence of a mutualist. The results of this section can be summarized in the following theorem.
\begin{theorem}\label{T:invader} Suppose $\alpha,\beta,\gamma,a,m,l,b,L>0$ and $a<1+mL$.
Then the equilibrium $E_4$ of the system \eqref{CCM} is invadable by the first species.\\
$\bullet$ For $b<1$ the outcome of the invasion is the $3$-species coexistence $E_7^+$. The system has a single spreading speed equal to the wave speed of a traveling wave with the profile connecting $E_4$ to $E_7^+$.\\
$\bullet$ For $b\geq1$ the outcome of the invasion is the resident extinction $E_6$.\\
$\bullet$ In all cases the slowest spreading speed is at least $c_1=2\sqrt{\alpha D_1\Bigl(1-\frac{a}{1+mL}\Bigr)}$.
If in addition $D_2\leq2D_1$,\hspace{1em} $D_3<2D_1+\frac{\gamma}{\alpha}\left(1-\frac{a}{1+mL}\right)^{-1}\!D_1$\hspace{1em} and
\begin{equation}\label{fsecineq}
\frac{ab-(1+mL)}{(1+mL)-a}\leq
\frac{\alpha}{\beta}\left(2-\frac{D_2}{D_1}\right)-\frac{\gamma}{\beta}\frac{mal}{(1+mL)^2}
\frac{\alpha\left(2-\frac{D_2}{D_1}\right)\left(1-\frac{a}{1+mL}\right)+\beta}
{\alpha\left(2-\frac{D_3}{D_1}\right)\left(1-\frac{a}{1+mL}\right)+\gamma}
\end{equation}
then the system is linearly determinate with a single spreading speed $c_1$.
\end{theorem}
\begin{proof}
The first two bullets follow from the above discussion and \refL{a1mL}. It remains to verify the conditions of \refT{LDet}.
By the above discussion, $a<1+mL$ implies condition I) of the theorem, so $c^*\geq c_1$ in all cases. For II) we need:
\begin{align}\label{gammadif}
\gamma_1(\overline{\mu})-\gamma_2(\overline{\mu})=
\alpha\left(2-\frac{D_2}{D_1}\right)\left(1-\frac{a}{1+mL}\right)+\beta>0\,; \notag\\
\gamma_1(\overline{\mu})-\gamma_3(\overline{\mu})=
\alpha\left(2-\frac{D_3}{D_1}\right)\left(1-\frac{a}{1+mL}\right)+\gamma>0\,.
\end{align}
The first inequality will follow from condition III) in our case, see below, and the second one gives the claimed inequality for $D_3$.
To check III) we need a technical result proved in the Appendix.
\begin{lemma}\label{L:g<g'} Let $g(q)$ be given by \eqref{gCCMI} and $\xi\gg0$ be arbitrary.
Then $g(\rho\xi)\leq\rho g'(0)\xi$ for all $\rho>0$ if and only if
\begin{equation}\label{xineq}
\begin{cases}b\xi_1\geq\xi_2\\
\xi_1\geq\frac{a\xi_2}{1+mL}+\frac{am\xi_3}{(1+mL)^2}\,.
\end{cases}
\end{equation}
\end{lemma}
\noindent We wish to make \eqref{xineq} explicit for $\xi=\zeta_1(\overline{\mu})$.
The first inequality from \eqref{xineq} reduces to $\gamma_1(\overline{\mu})-\gamma_2(\overline{\mu})\geq\beta$ or
$D_2\leq2D_1$ since we are already assuming $a<1+mL$. It also implies the first inequality from \eqref{gammadif} since $\beta>0$. The second inequality is more cumbersome:
\begin{equation}\label{secineq}
\gamma_1(\overline{\mu})-\gamma_2(\overline{\mu})\geq
\frac{\beta\,ab}{1+mL}+\frac{\gamma\,aml}{(1+mL)^2}\cdot
\frac{\gamma_1(\overline{\mu})-\gamma_2(\overline{\mu})}{\gamma_1(\overline{\mu})-\gamma_3(\overline{\mu})}\,.
\end{equation}
Substituting \eqref{gammadif} into \eqref{secineq} gives \eqref{fsecineq}.
\end{proof}
The first condition for linear determinacy, $D_2\leq2D_1$, is identical to the inequality for the $2$-species model,
see \cite{W2'}. The second one is vacuous for it (no mutualist), and \eqref{fsecineq} also reduces to the inequality for the $2$-species
model after setting $m=0$, as expected. For small $L$ at least, the range of linear determinacy seems to shrink as $m$ increases,
but keep in mind that \eqref{fsecineq} is only a sufficient condition.

\section{Mutualist-resident}\label{s5}

This section is structured identically to the previous one, while
replacing a mutualist-invader with a mutualist-resident. The results
however, are quite different. The linear spreading speed is not
affected by the mutualist, and the linear determinacy conditions of
\refT{LDet} are almost never met. We use a comparison principle
for cooperative systems to derive alternative estimates for the
spreading speeds.

Recall that for the mutualist-resident the source equilibrium is $E_6\,(1,0,L+l)$ with the invader absent. 
The converting substitution is
$$\ds{\left[\begin{array}{l} q_1=p_2\\
q_2=1-p_1\\
v=L+l-u\ (q_3=L+l-p_3),
\end{array}\right.}$$
where we switched the order to keep the invading species first. Correspondingly, we now set
$D:=\text{diag}(D_2,D_1,D_3)$ and obtain $\dot{q}=h(q)+Dq_{xx}$ in the vector form with
\begin{equation}\label{hCCMR}
h(q)=\begin{pmatrix}\beta\,q_1(1-b-q_1+bq_2)\\
\alpha\,(1-q_2)\left(\frac{aq_1}{1+mL+ml-mv}-q_2\right)\\
\gamma\,(L+l-v)\left(\frac{L+l-v}{L+l-lq_2}-1\right)\end{pmatrix}.
\end{equation}
The Jacobian matrix is $\ds{h'(q):=\left(\frac{\partial h_i}{\partial q_j}\right)}$
\begin{equation}\label{h'CCMR}
h'(q)=\begin{bmatrix}
\beta\,(1-b-2q_1+bq_2)& \beta\,b\,q_1 & 0\\
\alpha\,\frac{a(1-q_2)}{1+mL+ml-mv} & -\alpha\left(1-2q_2+\frac{aq_1}{1+mL+ml-mv}\right) &
\alpha\,\frac{maq_1(1-q_2)}{(1+mL+ml-mv)^2}\\
0 & \gamma\,l\frac{(L+l-v)^2}{(L+l-lq_2)^2} & \gamma\left(1-2\frac{L+l-v}{L+l-lq_2}\right)\end{bmatrix}.
\end{equation}
Again, for non-negative values of $\alpha,\beta,\gamma,a,b,m,L,l$ and $q_2\leq1$ the off-diagonal entries 
in \eqref{h'CCMR} are non-negative making the system cooperative.

The transformed equilibria are\\

\noindent\begin{tabular}{llll}\label{Gqbr}
$G_0\,(0,1,L+l)$ & $G_1\,(0,0,L+l)$ & $G_2\,(1,1,L+l)$ &  $G_3\,(0,1,l)$\\
$G_4\,(1,1,l)$ & $G_5\!\left(\frac{1-b}{1-ab},a\frac{1-b}{1-ab},L+l\right)$ & $G_6\,(0,0,0)$\,\, &  $G_7^\pm\,(1-bz_\pm,1-z_\pm,l(1-z_\pm))$\,.
\end{tabular}\\

\noindent The source equilibrium is of course the origin $G_6$. As far as the target equilibria are concerned, it is immediately clear that $G_7^-$ can not be one. Indeed, one can see by inspection that $G_7^+\ll G_7^-$ since $z_-<z_+$, so $G_7^-$ always has a strictly positive equilibrium between $0$ and itself. We will be able to say more about other possible targets in a moment.

For now, we compute $h'(0)$ from \eqref{h'CCMR}
\begin{equation}\label{h'0CCMR}
h'(0)=\begin{bmatrix}
\beta(1-b) & 0 & 0\\
\frac{\alpha\,a}{1+mL+ml}& -\alpha & 0\\
0 & \gamma\,l & -\gamma
\end{bmatrix}.
\end{equation}
This matrix is already in the Frobenius form with $1\times1$
diagonal blocks. Recall that we transposed the first and the
second species to achieve this. For positive $\alpha,\beta$, and
$\gamma$ a positive eigenvalue exists if and only if $b<1$, and we
assume $b<1$ throughout this section. The following lemmas show how this assumption influences the existence and 
the relative location of the equilibria. Let $p_i(E_j)$ denote the population density of the $i^{\text{th}}$ species in the equilibrium $E_j$.
\begin{lemma}\label{L:b<1-E5} Suppose $a,m,l,b,L>0$ and $b<1$. If $E_5$ is admissible
then $a\leq1$, and also $E_7^+$ is admissible while $E_7^-$ is inadmissible. Moreover, $p_1(E_5)<p_1(E_7^+)$,\quad $p_2(E_5)>p_2(E_7^+)$, in other words, the presence of a mutualist increases $p_1$ and decreases $p_2$ in the coexistence equilibrium.
\end{lemma}
\noindent Denote $B:=1+mL-ml-ab$ and $\mathcal{D}:=B^2+4ml(1+mL-a)$, then solutions to \eqref{zEq} are $z_{\pm}:=(-B\pm\sqrt{\mathcal{D}})/2ml$ and $z_+\leq z_-$ for $\mathcal{D}\geq0$. In the next lemma we treat the discriminant $\mathcal{D}$ as a quadratic polynomial in $a$.
\begin{lemma}\label{L:b<1} Suppose $a,m,l,b,L>0$ and $b<1$, $a>1+mL$. Then
$E_7^\pm$ are either both admissible or both inadmissible. Specifically, $\mathcal{D}(a)=0$ has two
positive real roots $a_1<a_2$, and $E_7^\pm$ are admissible if and only if
$a\leq a_1$ and $ml>(1+mL)(1-b)$. If $a=a_1$ then $E_7^+=E_7^-$, otherwise they are distinct.
\end{lemma}
\begin{figure}[!ht]\label{1}
\begin{centering}
\ \includegraphics[scale=0.9]{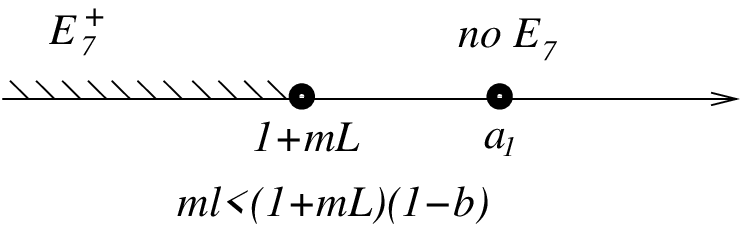} \hspace{0.2in} \ \includegraphics[scale=0.9]{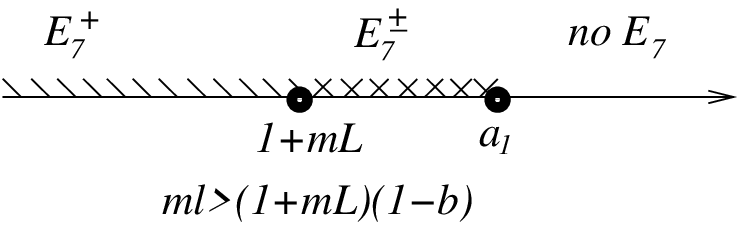}
\caption{\label{F:E7diagram} The existence of equilibria $E_7^\pm$ for $0<b<1$}
\end{centering}
\end{figure}
\noindent Lemmas \ref{L:b<1-E5} and \ref{L:b<1} refine the analysis in \cite{RFA}, their proofs are given in the Appendix.
The results of Lemmas \ref{L:a1mL} and \ref{L:b<1} are summarized on \refF{E7diagram}. We see that for $0<b<1$ there is always a change in the status of $E_7^\pm$ when $a$ crosses the value $1+mL$. Below it only $E_7^+$ is admissible, and above it either $E_7^+$ ceases to be admissible or $E_7^-$ also becomes admissible. Both equilibria $E_7^\pm$ cease to be admissible when $a$ crosses the larger value $a_1$. According to \cite{RFA}, in the conditions of \refL{b<1} the equilibrium $E_7^-$ is always unstable.

We are now ready to say more about the target equilibria in the mutualist-resident case. 
First, only $G_4$ and $G_7^+$ can
strictly dominate $G_6$ with no other equilibria in between.
Indeed, $G_7^+\ll G_5$ for $a\leq1$ by \refL{b<1-E5}, and for $a>1$ the
equilibrium $G_5$ is not admissible. We also have $G_7^+\ll G_4$ for
$a\leq1+mL$ by \refL{a1mL}, so $G_4$ can only be the target for
$a>1+mL$. By \refL{b<1}, it is definitely the target for $a>a_1$, and between $1+mL$ and $a_1$ either $G_4$ or $G_7^+$ is the target, depending on the sign of $ml-(1+mL)(1-b)$.
In all cases $q_2\leq1$ in the target equilibrium ensuring that the system remains cooperative throughout the evolution.

In the original system, $E_4\,(0,1,L)$ is the resident extinction
and  $E_7^+\,(z_+,1-bz_+,L+lz_+)$ is again the $3$-species coexistence. We
see that the extinction only happens for high competition intensities
$a$. Also note that the invader can not drive the mutualist to
extinction in this model due to a lack of direct interaction. If
$E_7$ is the target then the system has no boundary equilibria.
Hence, there is a single spreading speed equal to the wave speed
of a profile connecting $E_6$ to $E_7^+$. If $E_4$ is the target
then a boundary equilibrium $E_3$ is acquired, and no conclusion on a single speed follows.

We now proceed with the linear determinacy analysis. The matrix $C_\mu:=h'(0)+\mu^2D$ from \eqref{Cmu} has the following
diagonal entries, which are also its eigenvalues (recall that we transposed the first and the second
species in the conversion):
$$
\gamma_1(\mu)=\beta(1-b)+D_2\mu^2;\hspace{2em} \gamma_2(\mu)=-\alpha+D_1\mu^2;\hspace{2em}
\gamma_3(\mu)=-\gamma+D_3\mu^2.
$$
The principal eigenvalue is $\gamma_1$ with the eigenvector
\begin{equation}\label{z1(0)}
\zeta_1=\begin{pmatrix}
(\gamma_1-\gamma_2)(\gamma_1-\gamma_3)\\
\frac{\alpha\,a}{1+mL+ml}\,(\gamma_1-\gamma_3)\\
\frac{\alpha\,a}{1+mL+ml}\,\gamma\,l
\end{pmatrix},
\end{equation}
where we suppressed the dependence on $\mu$. By inspection, $\zeta_1(0)$ has positive entries if $b<1$.

Assuming the system is linearly determinate the linear speed is
\begin{equation}\label{c1r}
\ds{c_1=2\sqrt{\beta D_2(1-b)}}
\end{equation}
with
\begin{equation}\label{Rmubar}
\overline{\mu}=\sqrt{\frac{\beta}{D_2}\left(1-b\right)}.
\end{equation}
In contrast to the case of a mutualist-invader, this speed is the
same as for the $2$-species competition model of \cite{W2'} with
no mutualist present. An intuitive explanation is that in our case
the mutualist does not directly interact with the invader. Its
influence is indirect, it reduces the competition of the resident
with the invader, and is not captured by the linearized speed. In
another contrast, the linear determinacy condition III) of \refT{LDet} is
almost never satisfied as we will see from the next lemma.
\begin{lemma}\label{L:h<h'} Let $h(q)$ be given by \eqref{hCCMR} and $\xi\gg0$ be arbitrary. Suppose
$\alpha,\beta,\gamma,a,m,l>0$ and $0<b<1$. Then the third component of the vector inequality
$h(\rho\xi)\leq\rho h'(0)\xi$ holds for all $\rho>0$ and $\xi\gg0$ if and only if $\xi_3=l\xi_2$.
\end{lemma}
\noindent For $\xi=\zeta_1(\overline{\mu})$ the equality $\xi_3=l\xi_2$ becomes
$$
\frac{\alpha\,a}{1+mL+ml}\,\gamma\,l=\frac{\alpha\,a}{1+mL+ml}\,l\,
\left(\,\gamma_1(\overline{\mu})-\gamma_3(\overline{\mu})\,\right).
$$
Canceling $\alpha,a,l$ and taking into account \eqref{Rmubar}
turns this into
$$
\gamma=\gamma_1(\overline{\mu})-\gamma_3(\overline{\mu})
=\beta(1-b)\left(2-\frac{D_3}{D_2}\right)+\gamma\,,
$$
which is only possible if $D_3=2D_2$. Realistically, such precise equality between the 
mobilities of an invader and a mutualist almost never occurs. For this reason we do not derive
the inequalities for the other two components, which do lead to reasonable constraints.
Again, this is in stark contrast with the mutualist-invader case, where the inequality for 
the third component is always satisfied.

Thus, we have to take a different approach to estimating spreading
speeds here. Note that condition I) alone guarantees that the the
slowest speed $\ds{c^*\geq c_1=2\sqrt{\beta D_2(1-b)}}$. To get an
estimate from above we will compare our system
$\dot{q}=h(q)+Dq_{xx}$ to $\dot{q}=h^\circ(q)+Dq_{xx}$, where
$h^\circ(q)$ is obtained from $h(q)$ by setting $m=0$, i.e.
turning off the interaction with the mutualist. In terms of the
original model \eqref{CCM}, the first two equations then decouple
and form the usual diffusive Lotka-Volterra $2$-species
competition system \eqref{2LotVl}. Under the linear determinacy
conditions for \eqref{2LotVl} the spreading speed is again $c_1$ from \eqref{c1r}
giving an estimate from above for the original one. This means
that the invader spreads, and the resident recedes, at the same
speed $c_1$, leaving the possibility that the mutualist recedes
even faster. If we know in advance that the system does have a
single spreading speed, e.g. if $E_7^+$ is the target equilibrium,
then it will also be linearly determinate. Let us summarize our
discussion in a theorem.
\begin{theorem}\label{T:resident} Suppose $\alpha,\beta,\gamma,a,m,l,b,L>0$ and $b<1$.
Then the equilibrium $E_6$ of system \eqref{CCM} is invadable by the second species.\\
$\bullet$ Suppose $ml\leq(1+mL)(1-b)$ and $a\leq1+mL$, or $ml>(1+mL)(1-b)$ and $a\leq a_1$ with
$$
a_1=\frac1{b^2}\left(b(1+mL)+(2-b)ml-2\sqrt{ml(1-b)\left(b(1+mL)+ml\right)}\right)\,.
$$
Then the outcome of the invasion is $3$-species coexistence $E_7^+$, and there is a single spreading speed equal to the wave speed of the profile connecting $E_6$ to $E_7^+$.\\
$\bullet$ Suppose $ml\leq(1+mL)(1-b)$ and $a>1+mL$, or $ml>(1+mL)(1-b)$ and $a>a_1$.
Then the outcome of the invasion is the resident extinction $E_6$.\\
$\bullet$ In all cases the slowest spreading speed is at least $c_1=2\sqrt{\beta D_2(1-b)}$.
If in addition $D_1\leq2D_2$ and
\begin{equation}\label{hsecineq}
\frac{ab-1}{1-b}\leq
\frac{\beta}{\alpha}\left(2-\frac{D_1}{D_2}\right)
\end{equation}
then both the invader and the resident, but not necessarily the mutualist, spread (recede) at this speed. If moreover,
the outcome is $E_7^+$ then the system is linearly determinate with a single spreading speed $c_1$.
\end{theorem}
\begin{proof}
First two bullets follow from the above discussion and Lemmas \ref{L:a1mL}, \ref{L:b<1}.
Also by the above discussion, $b<1$ implies condition I) of \refT{LDet}, so $c^*\geq c_1$ in all cases.
Define $h^\circ(q)$ by setting $m=0$ in \eqref{hCCMR}:
\begin{equation}\label{h0CCMR}
h^\circ(q)=\begin{pmatrix}\beta\,q_1(1-b-q_1+bq_2)\\
\alpha\,(1-q_2)\left(aq_1-q_2\right)\\
\gamma\,(L+l-v)\left(\frac{L+l-v}{L+l-lq_2}-1\right)\end{pmatrix}\,.
\end{equation}
By inspection, $h(q)\leq h^\circ(q)$ for $0\leq
q\leq\beta\leq(1,1,l)$. Let $q(t,x)$, $q^\circ(t,x)$ be the
solutions to $\dot{q}=h(q)+Dq_{xx}$, $\dot{q}=h^\circ(q)+Dq_{xx}$
respectively with $q^\circ(0,x)=q(0,x)$. Then by the comparison
principle for cooperative systems \cite{Sm}, we have $q(t,x)\leq
q^\circ(t,x)$ for $t\geq0$. As remarked in \refS{s3}, this means
that the spreading speeds in the original system are bounded from
above by those in the $m=0$ system. But in the latter the first
two equations decouple and correspond to the diffusive
Lotka-Volterra $2$-species competition system (up to the
cooperativity substitution). Using the results of \cite{W2'} for
such systems, or setting $m=0$ and swapping $\alpha$ with $\beta$,
$a$ with $b$, and $D_1$ with $D_2$ in \refT{invader}, we get the
stated inequalities and the linear speed $c_1$ for the $2$-species
model. Since it coincides with the lower estimate above we
conclude that the invader and the resident spread (recede) at this
speed. From the results of \cite{W5} we also know that the
original system has a single speed if there are no boundary
equilibria between $0$ and $\beta$. This is the case when $E_7^+$ is
the target equilibrium. Together with $c_1$ being the spreading
speed for the first two components this guarantees linear
determinacy with a single speed $c_1$.
\end{proof}

Note that the comparison trick does not work in the mutualist-invader case since the inequality between $g$ and $g^\circ$
points in the wrong direction. It is an interesting question whether there actually are situations when the mutualist has
the spreading speed larger than $c_1$, while the invader and the resident spread at this speed. The failure of condition III) from \refT{LDet} in the mutualist's component only is certainly suggestive.

\section{Numerical simulations}\label{s6}

In this section we present the results of numerical simulations performed to measure the spreading speeds in 
different situations. 
\begin{figure}[!ht]
\begin{centering}
\includegraphics[scale=0.3]{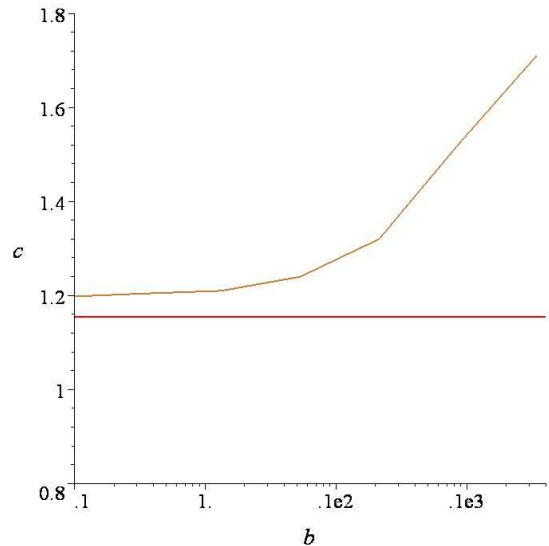}
\par\end{centering}
\caption{\label{figure1}The spreading block speeds $c_1(b)$ of the invader as a function of $b$ in the mutualist-invader case, see \refT{invader}. Here $a=\frac{2}{3}$ and $m=0$, $b$ varies on a log scale. The red line is the theoretical linear speed \eqref{c1i}, the gold line is a numerical estimate of the actual speed.}
\end{figure}
To estimate them we solved initial-boundary problems for the system \eqref{CCM} numerically on a long interval with a software that uses the method of lines. The initial profiles of the resident's and the mutualist's population densities were constant at the source equilibrium values, the boundary conditions were set to the source equilibrium values as well. The initial profile of the invader density was chosen as a narrow Gaussian bell curve centered at the middle of the interval. The spreading speeds were estimated by comparing spatial profiles of the densities at different times. The times chosen were large enough to allow for the speeds to settle, but not so large that boundary effects at the ends of the interval become substantial.

In all simulations below the values of the following parameters were fixed at $l=\frac{9}{20}, L=\frac{1}{3}, \alpha=1, \beta=1, \gamma=1, D_1=D_2=D_3=1$. Recall that the coefficients $a$ and $b$ in \eqref{CCM} are the intensities of competition between the invading and the resident species. For the mutualist-invader case Fig.\ref{figure1} shows the dependence of the spreading speed of the invader on the value of $b$ with $a=\frac{2}{3}$ and $m=0$, i.e. when the mutualist has no effect on the competition between the invader and the resident. This is essentially the Lotka-Volterra 2-species competition model with diffusion. 
\begin{figure}[!ht]
\begin{centering}
{\bf (a)} \includegraphics[scale=0.3]{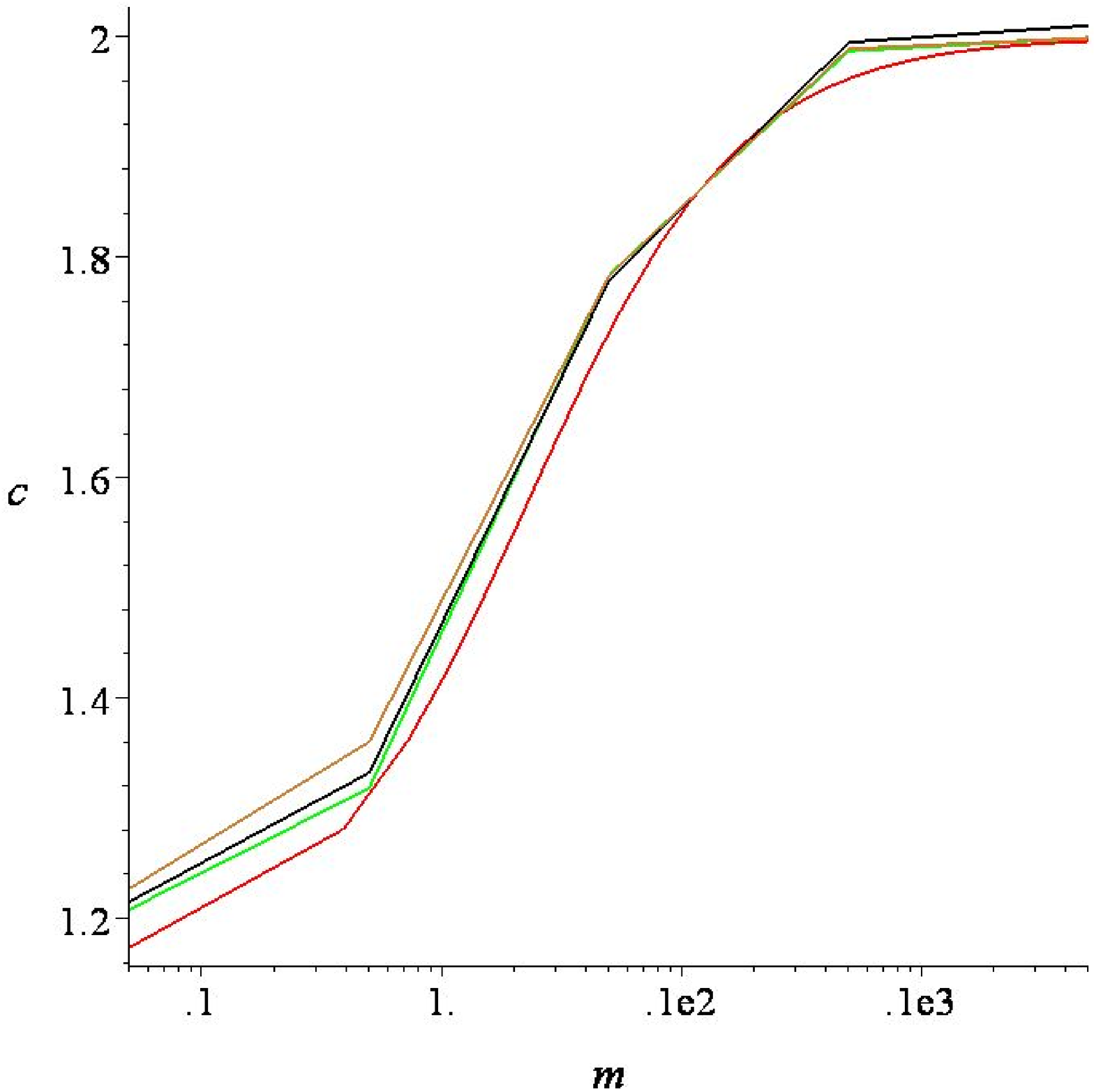}{\bf (b)} \includegraphics[scale=0.3]{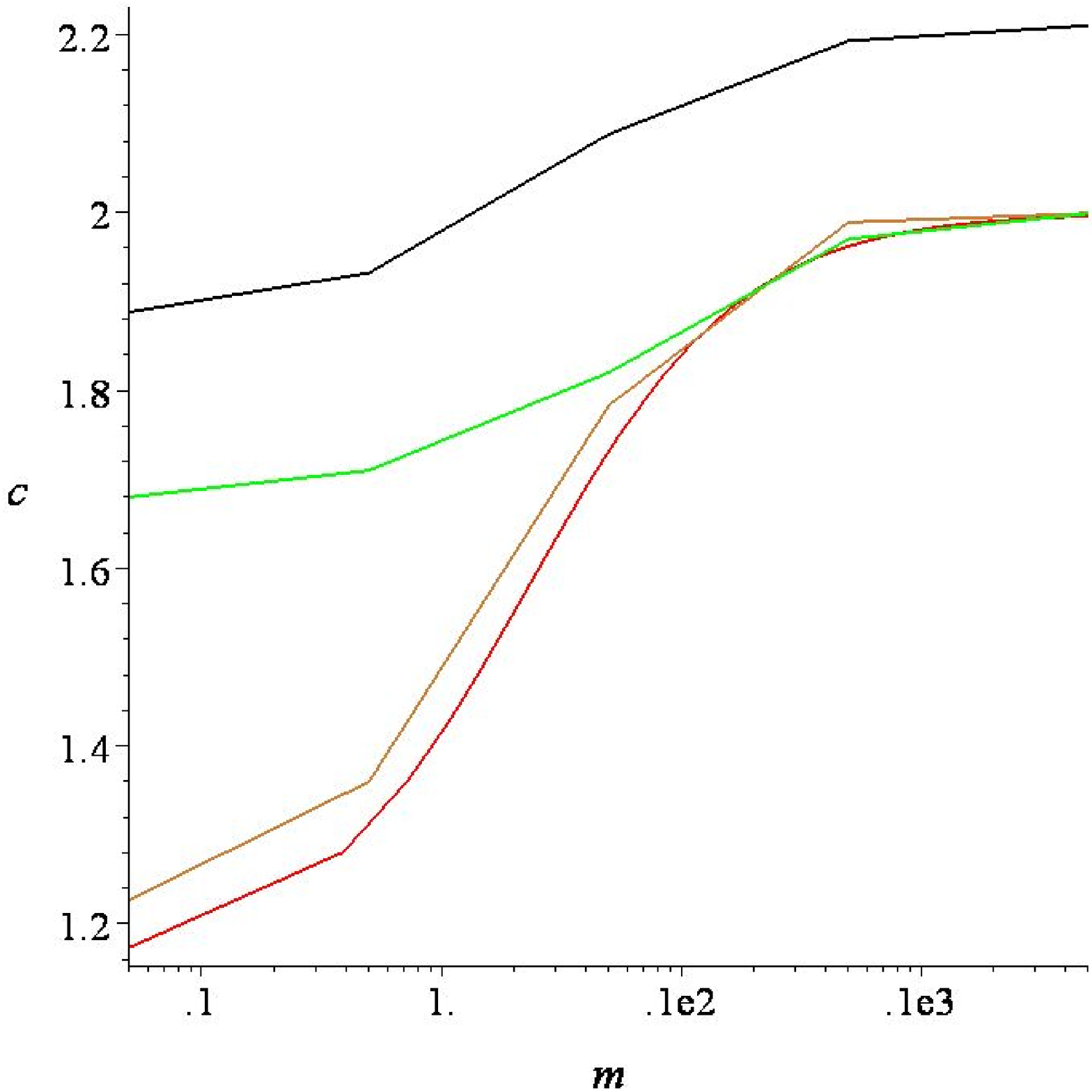}
\end{centering}
\caption{\label{figure2}Spreading block speeds $c_i(m)$ as functions of $m$ in the mutualist-invader case 
(\refT{invader}) for {\bf (a)} $b=\frac{2}{3}$ and {\bf (b)} $b=\frac{1024}{3}$ with $a=\frac{2}{3}$ in both cases, $m$ varies on a log scale. The red curve is the theoretical linear speed \eqref{c1i}, the rest are numerical estimates of the actual speeds: the invader's (gold), the resident's (black) and the mutualist's (green).}
\end{figure}
One can see that within a margin of error the estimated speed coincides with the linear speed for small $b$. With our values Theorem \ref{T:invader} guarantees linear determinacy for $b\leq1.5$, close to the value where the estimated speed is seen to start curving upward. This is consistent with the numerical results in \cite{W2'}. Note that at $b=1$ the target equilibrium switches from $E_7^+$ to $E_6$, but the linear speed does not change.

On Fig.\ref{figure2} we plotted the linear speed and the estimated spreading speeds of each species for two different values of $b$ as the intensity of mutualism $m$ increases. The smaller value $b=\frac{2}{3}$ keeps the system within the range of linear determinacy given by Theorem \ref{T:invader}, and one can see that the estimated speeds do indeed conform to the linear one. For the larger value $b=\frac{1024}{3}$ the invader's and the mutualist's speed appear to converge to the linear value asymptotically, but the speed of the resident's extinction remains consistently faster.
\begin{figure}[!ht]
\begin{centering}
{\bf (a)} \includegraphics[scale=0.3]{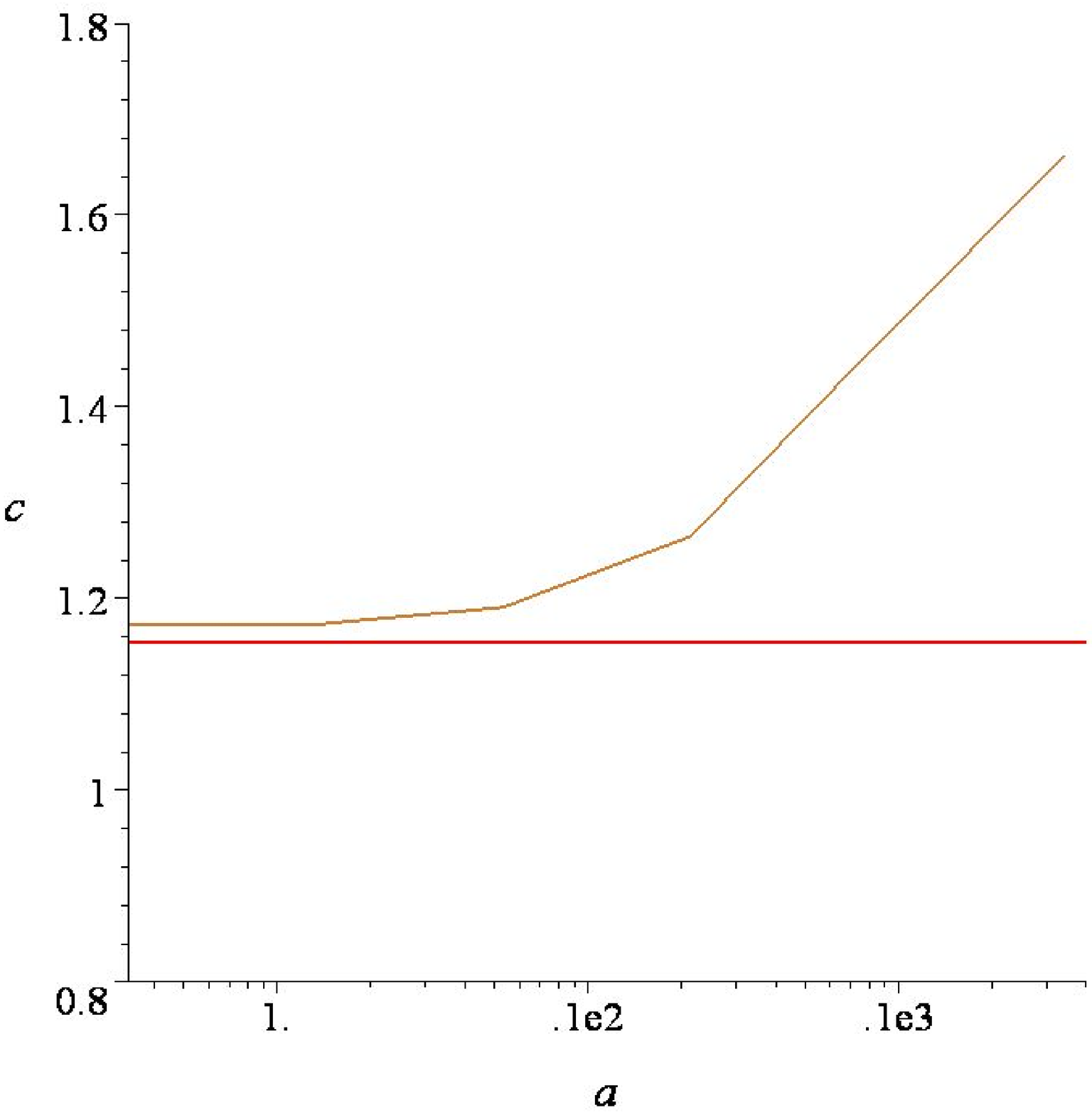}{\bf (b)} \includegraphics[scale=0.3]{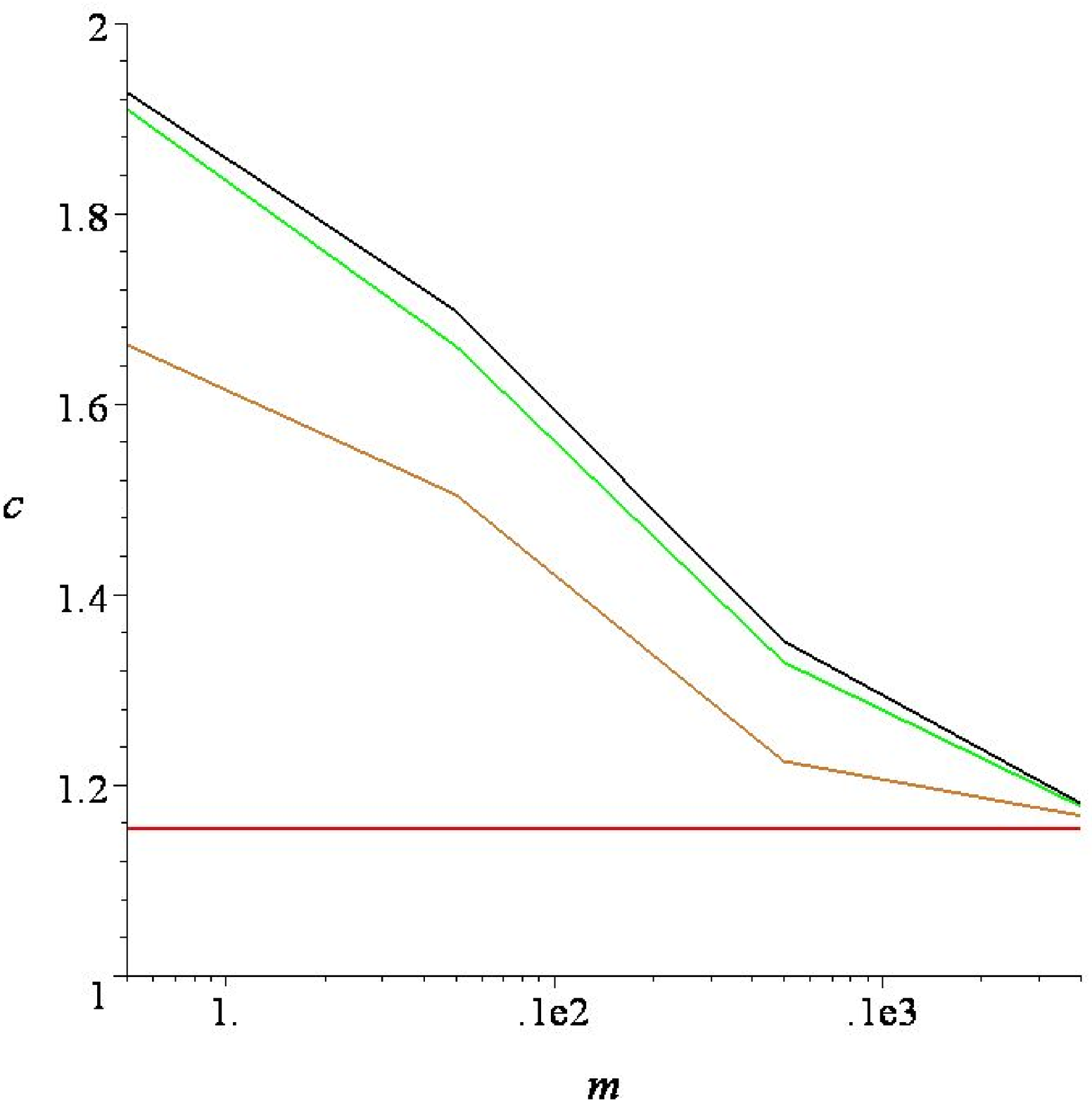}
\par\end{centering}
\caption{\label{figure3}{\bf (a)} The spreading block speeds $c_1(a)$ of the invader as functions of $a$ in the mutualist-resident case, see \refT{resident}. Here $b=\frac{2}{3}, m=\frac{1}{2}$ and $a$ varies on a log scale. The red line is the theoretical linear speed \eqref{c1r}, the gold line is the a numerical estimate of the actual speed. 
{\bf (b)} The spreading block speeds $c_i(m)$ as functions of $m$ in the mutualist-resident case for $b=2/3$ and $a=\frac{1024}{3}$, $m$ varies on a log scale. The red curve is the theoretical linear speed \eqref{c1r}, the rest are numerical estimates of the actual speeds: the invader's (gold), the resident's (black) and the mutualist's (green).}
\end{figure}

For the mutualist-resident case, Fig.\ref{figure3}(a) shows the invader's spreading speed as a function
of $a$ with $b=\frac{2}{3}$ and $m=\frac{1}{2}$. The behavior is similar to the situation on Fig.\ref{figure1} for the mutualist-invader: the speed retains the linear value for small $a$ and curves upward as $a$ increases. Theorem 
\ref{T:resident} guarantees linear determinacy for $a\leq1.5$. On Fig.\ref{figure3}(b) the value of $a$ is chosen to be beyond the range of linear determinacy. The speeds of the resident and the mutualist track each other closely, but differ substantially from the invader's speed. All speeds remain above the linear value, as the theory predicts, but appear to approach it from above asymptotically.

\section{Discussion}\label{s7}

We described two types of the invasion waves occurring in the
diffusive Competitor-Competitor-Mutualist system, in one the
mutualist benefits the invading, and in the other it benefits the
resident species. In each case we obtained estimates for the
spreading speeds and sufficient conditions for linear determinacy.
In both cases the mutualist alters the coexistence outcome in favor of the species it benefits. In
the mutualist-invader case the speed of the invasion (under linear
determinacy) is increased  compared to the diffusive
Lotka-Volterra model with no mutualist. In the mutualist-resident
case this speed remains the same due to lack of direct interaction
between the invader and the mutualist in the model. While for the
mutualist-invader the linear determinacy conditions are direct
generalizations of those for the Lotka-Volterra model, for the
mutualist-resident such direct generalization fails, and the
nature of the spreading speed is more subtle. In both cases these
conditions can be interpreted as requiring relatively high
mobility of the invader compared to the resident (and in the first
case, to the mutualist), and relatively weak competition between
the invader and the resident.

It is instructive to collect different invasion scenarios on a single diagram, \refF{outcomes}.
\begin{figure}[!ht]
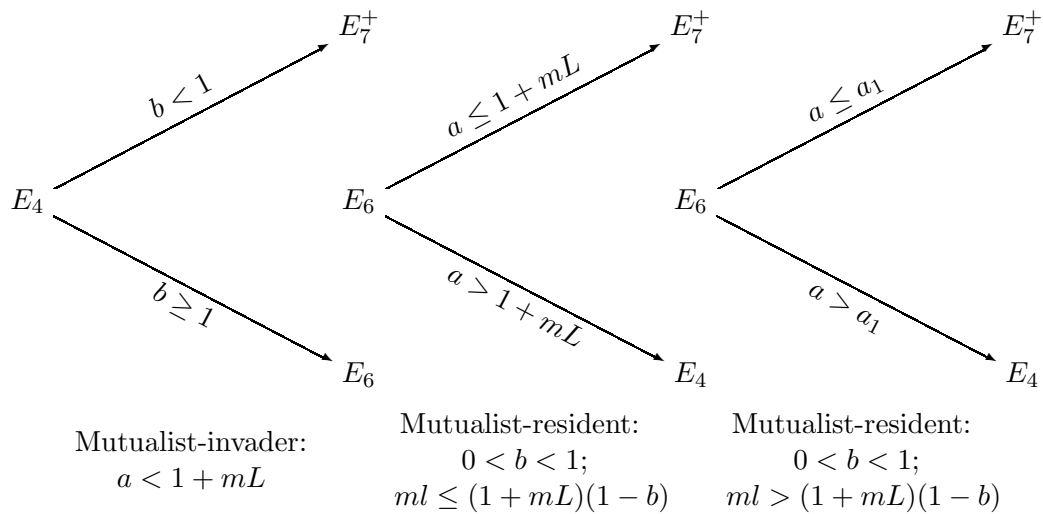

\begin{equation*}
\begin{diagram}
   &                & E_7^+ &      &           E_7^+& &                    E_7^+ \\
   & \ruTo^{\textstyle{b<1}}    &     &       \ruTo^{\textstyle{a\leq1+mL}}  &  &    \ruTo^{\textstyle{a\leq a_1}}&         \\
E_4&                &       E_6 & &                              E_6&  & &                     \\
   & \rdTo_{\textstyle{b\geq1}} &     &       \rdTo_{\textstyle{a>1+mL}}  &        &     \rdTo_{\textstyle{a>a_1}}&        \\
   &                & E_6 &      &                          E_4 & &                  E_4 \\
   &  \begin{matrix}\textstyle{\text{Mutualist-invader:}}\\\hspace{2.5em}\textstyle{a<1+mL} \hspace{2.5em} \end{matrix}      &     &       \begin{matrix}\textstyle{\text{Mutualist-resident:\ }}\\\textstyle{0<b<1};\\\,\,\,\textstyle{ml\leq(1+mL)(1-b)}\end{matrix}& &     \begin{matrix}\textstyle{\text{Mutualist-resident:\ }}\\\textstyle{0<b<1};\\
\,\,\,\textstyle{ml>(1+mL)(1-b)} \end{matrix}
\end{diagram}
\end{equation*}
\caption{\label{F:outcomes} Diagram of invasion outcomes}
\end{figure}
Comparing the results for the mutualist-invader and the
mutualist-resident we are led to an interesting observation. If
the system starts at $E_6$ and $a>1+mL$ then the first species
will be wiped out by the invasion, resulting in $E_4$. But $E_4$ is
also the source equilibrium for the invasion (or in this case,
reinvasion) by the first species discussed in \refS{s4}. Of
course, one needs $a<1+mL$ for the reinvasion to be possible,
contrary to the above. However, if one can increase the
self-carrying capacity of the mutualist to $L'$ so that $a<1+mL'$,
then the reinvasion wave becomes possible and leads to the
$3$-species coexistence $E_7^+$. This indicates that one can
successfully re-establish the native species by manipulating the
mutualist only, an attractive property if it is used as a
biocontrol agent.
\bigskip

{\em Acknowledgements:} The idea of this work was conceived when the first author attended the NIMBioS investigative workshop New Soil Black Box Math Strategies. He wishes to thank the organizers for their hospitality and the stimulating atmosphere, and Alan Hastings for introducing him to mathematical aspects of the spreading speeds. 

\section*{Appendix}

\begin{proof}[Proof of {\bf \refL{a1mL}}] Since $z$, $1-bz$ are the population densities in $E_7^\pm$ only non-negative values for them are admissible, so $0\leq z\leq1/b$. Recall that
\begin{equation}\label{BDz}
B:=1+mL-ml-ab;\hspace{.3em}\mathcal{D}:=B^2+4ml(1+mL-a);\hspace{.3em}z_{\pm}:=\frac{-B\pm\sqrt{\mathcal{D}}}{2ml}\,.
\end{equation}

By assumption, $1+mL-a>0$, so $\mathcal{D}>0$ and
$|B|<\sqrt{\mathcal{D}}$. Therefore, $z_+$ is the only non-negative solution.
Next, we show that $1-z_+$ and $1-b$ have the same sign:
\begin{multline}\label{1-z+}
1-z_+=\frac{B+2ml-\sqrt{\mathcal{D}}}{2ml}=\frac{(B+2ml)^2-\mathcal{D}}{2ml\,(B+2ml+\sqrt{\mathcal{D}})}\\
=\frac{4ml(-ab+a)}{2ml\,(B+2ml+\sqrt{\mathcal{D}})}=\frac{2a(1-b)}{B+\sqrt{\mathcal{D}}+2ml}.
\end{multline}
Since $B+\sqrt{\mathcal{D}}\geq0$ we conclude that $z_+$ and $b$ are on the same side of $1$. In particular, if $b>1$ then $z_+>1$, so $1-bz_+<0$ and $E_7^+$ can not be admissible. On the other hand, if $b\leq1$ then 
$z_+\leq1$, and $1-bz_+\leq0$ implies admissibility.

By inspection from \eqref{Eqbr}, $E_7^+$ can only merge with other equilibria if $z_+$ is $0$ or $1$. By \eqref{1-z+} the latter happens if and only if $b=1$, in which case $E_7^+=E_6$. To have $z_+=0$ one needs the free 
term $a-1-mL$ to vanish, which is excluded by the assumption $a<1+mL$.
\end{proof}

\begin{proof}[Proof of {\bf \refL{g<g'}}] First, we derive inequalities equivalent to $g(q)\leq g'(0)q$ for
arbitrary $q\gg0$. Utilizing \eqref{gCCMI} and \eqref{g'0CCMI} after some algebra one gets for each component
respectively:
$$
\begin{cases}q_1(1+mL)(1+mL+mq_3)\geq aq_2(1+mL)+amq_3\\
bq_1-q_2\geq0\\
\frac{(lq_1-q_3)^2}{L+lq_1}\geq0.
\end{cases}
$$
The last inequality is vacuous and the second one yields $\rho b\xi_1\geq\rho\xi_2$ after setting $q=\rho\xi$. This is satisfied for all $\rho>0$ if and only if $b\xi_1\geq\xi_2$. After setting $q=\rho\xi$ in the first inequality and canceling $\rho$ one has
$$
\xi_1(1+mL)(1+mL+\rho m\xi_3)\geq a\xi_2(1+mL)+am\xi_3.
$$
By passing to limit, this holds for any $\rho>0$ if and only if it holds for $\rho=0$, i.e.
$$
\xi_1(1+mL)^2\geq a\xi_2(1+mL)+am\xi_3.
$$
Division by $(1+mL)^2$ yields the desired claim.
\end{proof}

\begin{proof}[Proof of {\bf \refL{b<1-E5}}] Since $E_5$ is admissible its entries must be non-negative, meaning that
$1-a$, $1-b$ and $1-ab$ have the same sign. Thus, if $b<1$ then $a\leq1$ and \refL{a1mL} implies that $E_7^+$ is admissible, while $E_7^-$ is not.

We see that $p_1(E_7^+)=z_+$, $p_1(E_5)=\frac{1-a}{1-ab}$ and $p_2(E_7^+)=1-bz_+$, $p_2(E_5)=\frac{1-b}{1-ab}$. Thus,
it remains to show that $z_+>\frac{1-a}{1-ab}\text{ and }1-bz_+<\frac{1-b}{1-ab}$. The second inequality follows from the first one since $\frac{1-b}{1-ab}=1-b\,\frac{1-a}{1-ab}$. To prove the first one separate the terms in \eqref{zEq} with and without $m$ to get
$$
-m(1-z)(L+lz)+(1-ab)z-(1-a)=0.
$$
Now 'solve' for $z$ in the middle term:
$$
z=\frac{1-a}{1-ab}+m(L+lz)\frac{1-z}{1-ab}.
$$
It follows from the proof of \refL{a1mL} that $z_+<1$ for $b<1$, 
so $\frac{1-z_+}{1-ab}>0$. This implies the first inequality.
\end{proof}

\begin{proof}[Proof of {\bf \refL{b<1}}] Recall that for $E_7^\pm$ to be admissible both $z_\pm$ and $1-bz_\pm$ should be non-negative yielding $0\leq z_\pm\leq1/b$ as the admissibility criterion. The proof is rather cumbersome and we split it into two steps.

\underline{Step 1.} In this step we prove that $z_\pm$ are both admissible if one is, and it happens if and only if $\mathcal{D}\geq0$ and $-2ml/b<B<0$. Note that $\mathcal{D}\geq0$ is necessary and sufficient for real roots to exist at all, so we proceed to prove the double inequality for $B$ assuming that both roots are real. 

Let $G(z)=ml\,\!z^2+(1-ab+mL-ml)z-(1+mL-a)$ be the left hand side of \eqref{zEq}. Its graph is a parabola that opens up since $ml>0$. Since $G(0)=-(1+mL-a)>0$ both of its roots are on the same side of $0$. Similarly, since 
$G(1/b)=(1/b-1)(ml/b+1+mL)>0$ they are on the same side of $1/b$ as well. For the equilibria to be admissible we must have 
$z_\pm\in[0,1/b]$. Therefore, one of them is admissible if and only if both are. Since the vertex of the parabola is the midpoint of $[z_-,z_+]$ the roots belong to $[0,1/b]$ if and only if the vertex does, i.e. $-B/2ml\in[0,1/b]$. Note that the endpoints are excluded since $G(0)$ and $G(1/b)$ are strictly positive. But then the admissibility condition becomes 
$-2ml/b<B<0$ as claimed.

\underline{Step 2.} We now solve the system of three inequalities $\mathcal{D}\geq0$ and $-2ml/b<B<0$. The double inequality for $B$ can be written explicitly as
\begin{equation}\label{Db<B}
\frac{1+mL-ml}b<a<\frac{1+mL+ml(2/b-1)}b.
\end{equation}
On the other hand, we see from \eqref{BDz} that $\mathcal{D}=\mathcal{D}(a)$ is quadratic in $a$ with
the leading coefficient $b^2>0$. By direct computation,
$$
\mathcal{D}\Bigr(\frac{1+mL+ml(2/b-1)}b\Bigl)=-\frac{4ml}b\,(1-b)(1+mL+ml/b)<0\,.
$$
Thus, $\mathcal{D}$ has two real roots $a_{1,2}$ satisfying $a_1<(1+mL+ml(2/b-1))/b<a_2$\,, and 
$\mathcal{D}(a)\geq0$ if and only if $a\leq a_1$ or $a\geq a_2$\,.

In view of \eqref{Db<B} all three inequalities are satisfied if and only if 
\begin{equation}\label{Dba}
\frac{1+mL-ml}b<a\leq a_1\,.
\end{equation} 
Solutions to this double inequality exist if and only if $\frac{1+mL-ml}b<a_1$, i.e.
\begin{equation}\label{Sat}
\mathcal{D}\Bigr(\frac{1+mL-ml}b\Bigl)=4ml\Bigr(1+mL-\frac{1+mL-ml}b\Bigl)>0\,,
\end{equation}
which reduces to $1+mL>\frac{1+mL-ml}b$. Since $a>1+mL$ by assumption of the Lemma only $a\leq a_1$ remains to be required in \eqref{Dba} whenever \eqref{Sat} holds. But the latter is equivalent to $ml>(1+mL)(1-b)$.

The last claim is obvious since $z_\pm$ merge if and only if $\mathcal{D}(a)=0$, and $a=a_2$ keeps $z$ always outside the range of admissibility.
\end{proof}

\begin{proof}[Proof of {\bf \refL{h<h'}}] From \eqref{hCCMR} and \eqref{h'0CCMR} we have for the third component
$$
\left(\rho\,l\xi_2-\rho\xi_3\right)\,\frac{L+l-\rho\xi_3}{L+l-\rho\,l\xi_2}\left(\rho\,l\xi_2-\rho\xi_3\right)\leq
\gamma\left(\rho\,l\xi_2-\rho\xi_3\right).
$$
Canceling $\gamma$, $\rho$ and collecting all terms on the left:
$$
\left(l\xi_2-\xi_3\right)\left(\frac{L+l-\rho\xi_3}{L+l-\rho\,l\xi_2}-1\right)
=\rho\,\frac{\left(l\xi_2-\xi_3\right)^2}{L+l-\rho\,l\xi_2}\leq0.
$$
For small positive $\rho$ the denominator is positive, hence $l\xi_2-\xi_3=0$.
\end{proof}

\end{document}